\documentclass{amsart}

\newtheorem{theorem}{Theorem}[section]
\newtheorem{proposition}[theorem]{Proposition}
\newtheorem{corollary}[theorem]{Corollary}
\newtheorem{lemma}[theorem]{Lemma}

\newenvironment{nproof}[2]{\removelastskip\vspace{6pt}\noindent
 {\it Proof  #1.}~\rm#2}{\par\vspace{6pt}}

\theoremstyle{definition}

\theoremstyle{remark}
\newtheorem{remark}{Remark}

\numberwithin{equation}{section}

\DeclareMathOperator{\Tr}{Tr}

\begin{document}

\title[Linear relations among holomorphic quadratic differentials]{Linear relations among holomorphic quadratic \\
differentials and induced Siegel's metric on $\mathcal{M}_g$}


\author{Marco Matone}
\address{Dipartimento di Fisica ``G. Galilei'' and Istituto
Nazionale di Fisica Nucleare \\
Universit\`a di Padova, Via Marzolo, 8 --\\35131 Padova, Italy}
\email{matone@pd.infn.it} \email{volpato@pd.infn.it}
\thanks{Work partially supported by the European Community's Human Potential Programme under contract
MRTN-CT-2004-005104 ``Constituents, Fundamental Forces and
Symmetries of the Universe".}

\author{Roberto Volpato}


\date{September 2006}



\newcommand{\ZZ}{\mathbb{Z}}
\newcommand{\RR}{\mathbb{R}}
\newcommand{\CC}{\mathbb{C}}
\newcommand{\DD}{\mathbb{D}}
\newcommand{\II}{\mathbb{I}}
\newcommand{\QQ}{\mathbb{Q}}
\newcommand{\NN}{\mathbb{N}}
\newcommand{\KK}{\mathbb{K}}
\newcommand{\PP}{\mathbb{P}}
\newcommand{\FF}{\mathbb{F}}

\newcommand{\A}{\mathcal{A}}
\newcommand{\B}{\mathcal{B}}
\newcommand{\D}{\mathcal{D}}
\newcommand{\N}{\mathcal{N}}
\newcommand{\G}{\mathcal{G}}
\newcommand{\F}{\mathcal{F}}
\newcommand{\Sieg}{\mathfrak{H}}
\newcommand{\I}{\mathcal{I}}
\newcommand{\M}{\mathcal{M}}
\newcommand{\Z}{\mathcal{Z}}
\newcommand{\Os}{\mathcal{O}}
\newcommand{\1}{\mathfrak{1}}
\newcommand{\2}{\mathfrak{2}}
\newcommand{\perm}{\mathcal{P}}

\begin{abstract}
We derive the explicit form of the $(g-2)(g-3)/2$ linearly
independent relations among the products of pairs in a basis of
holomorphic abelian differentials in the case of canonical curves of
genus $g\geq4$. It turns out that Petri's relations
remarkably match in determinantal conditions. We explicitly express the
volume form on the moduli space ${\hat\M_g}$ of canonical curves induced by the
Siegel metric, in terms of the period Riemann matrix only.
By the Kodaira-Spencer map, the relations lead to an expression
of the induced Siegel metric on ${\hat\M_g}$, that corresponds
to the square of the Bergman reproducing kernel.
A key role is played by distinguished bases for holomorphic differentials whose properties
also lead to an immediate derivation of Fay's trisecant identity.

\end{abstract}
\maketitle

\section{Introduction}
\label{intro} In spite of the remarkable progresses in understanding
the Schottky problem, the characterization of the Schottky locus as
the zero set of modular forms on the Siegel's upper half-space
remains a fundamental open question. Such a question is strictly
related to the problem of characterizing the Schottky locus by means
of $(g-2)(g-3)/2$ linearly independent conditions. As suggested also
by Mumford (see pg. 241 of \cite{mumfordd}), a possible solution of
that problem should follow by a deeper understanding of Petri's
construction \cite{petriuno,ottimo}. Actually, since Petri's theorem
determines the ideal of canonical curves of genus $g\geq4$ by means
of linear relations among holomorphic abelian differentials, it
seems the natural framework for such an investigation.

Let $\{\eta_1,\ldots,\eta_g\}$ be the Petri's basis for $H^0(K_C)$,
with $C$ a canonical curve of genus $g$. In Petri's work the
coefficients $C_{ij}^k$ in the relationships among quadratic
differentials $\sum_{i,j}^g C_{ij}^k\eta_i\eta_j=0$, $k=1,\ldots,
(g-2)(g-3)/2$, are not determined. Finding such coefficients is a
necessary condition for an explicit characterization of the ideal of
canonical curves. Here, we express Petri's relation in determinantal
form, so that, besides the explicit determination of the
coefficients, it is shown that the locus of canonical curves
corresponds to a determinantal variety.

We introduce modular invariant bases for
holomorphic differentials, leading to a refinement of Petri's basis and to
an immediate derivation of Fay's trisecant identity \cite{Fay}.
A key point is the introduction of a indexing, which includes the combinatorics of the Petri
construction, mapping the components of matrices in the
Siegel upper half-space to vector components. This provides
the volume form on the moduli space ${\hat\M_g}$ of canonical curves induced by the
Siegel metric which, remarkably, is expressed in terms of the period
Riemann matrix only. By the Kodaira-Spencer map, the above relations lead to an expression
of the metric on ${\hat\M_g}$, induced by the Siegel metric, that corresponds
to the square of the Bergman reproducing kernel.

In the case of branched covering of the torus, corresponding to Jacobians
with a distinguished complex multiplication \cite{MatoneUY}, the derived relations should lead to identities of number theoretical
interest. Our results, of interest also in superstring theory 
\cite{MatoneVM}, provide the key for the $(g-2)(g-3)/2$
combinatorial $\theta$-identities in \cite{MatoneVO}.

\section{Determinantal characterization of canonical curves}
\label{sec:1}Let $C$ be a canonical curve of genus $g\ge 4$ and
$\{\omega_i\}_{i\in I_g}$, $I_n:=\{1,\ldots,n\}$, a basis of
$H^0(K_C)$, with $K_C$ the canonical line bundle of $C$. Denote
by $\hat\M_g$ the corresponding locus in the moduli space $\M_g$ of
compact Riemann surfaces. Each element of $H^0(K_C^2)$ can be
written as a linear combination of the $M:=g(g+1)/2$ elements in
$$
\mathcal{S}:=\{\omega_i\omega_j|i\leq j\in I_g\}\ .
$$
Since $N:= h^0(K_C^2)=3g-3$, there are $M-N=(g-2)(g-3)/2$ linearly
independent relations among the quadratic differentials
$\omega_i\omega_j$.

\noindent Let $p_1,\ldots,p_g$ and $q_1,\ldots,q_{2g-2}$ be two sets
of arbitrary points on $C$. Choose a local trivialization of the
canonical line bundle and set
$$a_{ij,r}:=\det\,\omega(p_1,\ldots,p_{i-1},q_r,p_{i+1},\ldots,p_g)\det\,\omega(p_1,\ldots,p_{j-1},q_r,p_{j+1},\ldots,p_g)\ ,$$
where $\det\,\omega(x_1,\ldots,x_g):= \det\, \omega_i(x_j)$. Set
$$A(k,l):=
\begin{pmatrix} a_{12,1} &\ldots &a_{1g,1} &a_{23,1} &\ldots &a_{2g,1} &a_{kl,1}   \\ a_{12,2} &\ldots &a_{1g,2} &a_{23,2}
&\ldots &a_{2g,2} &a_{kl,2}   \\ \vdots &\ddots  &\vdots
&\vdots &\ddots &\vdots &\vdots
\\ a_{12,2g-2} &\ldots
&a_{1g,2g-2}    &a_{23,2g-2} &\ldots &a_{2g,2g-2}
&a_{kl,2g-2}\end{pmatrix}\ ,$$ $3\le k<l\le g$, $g\ge 4$.

\begin{theorem}\label{teo}$$\det\,A(k,l)=0\ ,\qquad 3\le k<l\le g\ .$$\end{theorem}

\noindent Set $$\Delta_{m\, n}:=(-)^{m+n}\det_{^{i\neq m}_{j\neq n}}
A(k,l)_{ij}\ , \qquad D_{p\, q}:=(-)^{p+q}\det_{^{i\neq p}_{j\neq
q}}\omega_j(p_i)\ ,$$ and denote by $A_{ij,r}(k,l)$, $i,j\in I_g$,
$r\in I_{2g-2}$, the matrix obtained from $A(k,l)$ by replacing
the row $(a_{12,r},\ldots,a_{kl,r})$ with
$(D_{1i}D_{2j},\ldots,D_{k\, i} D_{l\, j})$.

\begin{corollary}\label{corolla} For each $r\in I_{2g-2}$, the following relations
$$\sum_{i,j=1}^g\frac{\det\,A_{ij,r}(k,l)}{\Delta_{r\, 2g-2}}\omega_i\omega_j=0\ ,$$
$3\le k<l\le g$, provide $(g-2)(g-3)/2$ linearly independent
conditions on $\mathcal{S}$ which are independent of the points
$q_i$, $i\in I_{2g-2}$. \end{corollary}

\section{Distinguished bases of $H^0(K_C^n)$}\label{sec:2}
Set $N_n:=(2n-1)(g-1)+\delta_{1n}$,
$n\ge 1$, with $\delta_{ij}$ the Kronecker delta. Note that
$N_1=g$ and $N_2\equiv N$. Fix a system of local
coordinates on $C$.

\begin{proposition}\label{prop:base} Fix $n\ge 1$ and let $p_1,\ldots,p_{N_n}$ be a set of points
in $C$ such that $$\det\,\phi(p_1,\ldots,p_{N_n})\ne 0\ ,$$ for
an arbitrary basis $\{\phi_i\}_{i\in I_{N_n}}$ of $H^0(K_C^n)$.
Then
\begin{equation}\label{basendiff}\gamma^n_i(z):=
\frac{\det\,\phi(p_1,\ldots,p_{i-1},z,p_{i+1},\ldots,p_{N_n})}{
\det\,\phi(p_1,\ldots,p_{N_n})}\ ,\end{equation}
$i\in I_{N_n}$, determines a basis of $H^0(K_C^n)$ which is
independent of the choice of the basis $\{\phi_i\}_{i\in
I_{N_n}}$ and, up to normalization, of the local coordinates on $C$.
\end{proposition}

\begin{proof} The matrix $[\phi]_{ij}:=\phi_i(p_j)$ is non-singular. Then
$\gamma^n_i(z)=\sum_j[\phi]^{-1}_{ij}\phi_j(z)$,
$i\in I_{N_n}$, is a basis of $H^0(K_C^n)$.\end{proof}

\noindent Note that $\gamma_i^n(p_j)=\delta_{ij}$, furthermore
\begin{equation}\label{preludetofay}
\det\, \gamma^n(p_1,\ldots,p_{j-1},z,p_{j+1},\ldots,p_{N_n})=\gamma^n_j(z) \ .
\end{equation}

\begin{remark} As we will see, the Fay trisecant identity \cite{Fay}
directly follows by expressing Eq.\eqref{basendiff} in terms of
theta functions.\end{remark}

For $n=1$, the choice of $g$ points $p_1,\ldots,p_g\in C$, with
$\det\,\omega_i(p_j)\ne 0$, determines the basis $\{\sigma_i\}_{i\in
I_g}$ of $H^0(K_C)$, where
\begin{equation}\label{newbasis}\sigma_i(z):=\gamma^1_i(z)\ ,\qquad i\in I_g \ .\end{equation}

We now introduce a refinement of Petri's basis for $H^0(K_C^2)$
\cite{petriuno,ottimo} which provides a modular
invariant construction. Let us assume that the points
$p_1,\ldots,p_g$ are in ``general position" and that the effective
divisor $(\sigma_1)+(\sigma_2)-\sum_{i=3}^gp_i$ consists of $3g-3$
distinct points. Consider the following $M$ elements of $H^0(K_C^2)$
$$v_i:=\begin{cases}\sigma_i^2\ , &i\in I_g\
,\\ &\\ \sigma_{j+k}\sigma_j\ , &i=k+j(2g-j+1)/2\ \ ,\quad
j\in I_{g-1}\ ,\quad k\in I_{g-j}\ . \end{cases}$$

\begin{proposition} $\{v_i\}_{i\in I_N}$ is a basis of $H^0(K_C^2)$.\end{proposition}

\begin{proof} Set $D:=\sum_{i=3}^{g}p_i$ and let us consider the effective divisors $D_i:=(\sigma_i)-D$, $i=1,2$. Let
us first prove that $\sigma_i$ is the unique $1$-differential, up to
normalization, vanishing at $D_i$, $i=1,2$. Any $1$-differential
$\sigma'_i\in H^0(K_C)$ vanishing at $D_i$, corresponds to an
element $\sigma'_i/\sigma_i$ of $H^0(\Os(D))$, the space of
meromorphic functions $f$ on $C$ such that $(f)+D$ is an effective
divisor. Suppose that there exists $\sigma'_i$ such that
$\sigma'_i/\sigma_i$ is not a constant, so that $h^0(\Os(D))\ge 2$.
By Riemann-Roch theorem $$h^0(K_C\otimes \Os(-D))=h^0(\Os(D))-\deg
D-1+g\ge 3\ ,$$ so that there exist at least $3$ linearly
independent $1$-differentials vanishing at $D$ and, in particular,
there exists a linear combination of such differentials vanishing at
$p_1,\ldots,p_g$. This implies that $\det\eta(p_1,\ldots,p_g)=0$ for
an arbitrary basis $\{\eta_i\}_{i\in I_g}$ of $H^0(K_C)$,
contradicting the hypotheses. Fix $\zeta_i,\zeta_{1i},\zeta_{2i}\in
\CC$ such that
$$\sum_{i=3}^g\zeta_i\sigma_i^2+\sum_{i=1}^g\zeta_{2i}
\sigma_1\sigma_i+\sum_{i=2}^g\zeta_{1i}\sigma_2\sigma_i=0\ .$$
Evaluating this relation at the point $p_j$, $3\le j\le g$, yields
$\zeta_j=0$. Set
\begin{equation}\label{leti}t_1:=-\sum_{j=2}^g\zeta_{1j}\sigma_j\ ,\quad
t_2:=\sum_{j=1}^g\zeta_{2j}\sigma_j\ ,\end{equation} so that
$\sigma_1t_2=\sigma_2 t_1$. Since $D$, $D_1$ and $D_2$ consist of
pairwise distinct points, $t_i$ vanishes at $D_i$, $i=1,2$ and then
$t_1/\sigma_1=t_2/\sigma_2=\zeta\in\CC$. By
\eqref{leti}
$$\zeta\sigma_1+\sum_{j=2}^g\zeta_{1j}\sigma_j=0\ ,\quad
\zeta\sigma_2-\sum_{k=1}^g\zeta_{2k}\sigma_k=0\ ,$$ and, by linear
independence of $\sigma_1,\ldots,\sigma_g$, we have
$\zeta=\zeta_{1j}=\zeta_{2k}=0$, $2\le j\le g$, $k\in I_g$.\end{proof}

\section{Proofs of Theorem \ref{teo} and Corollary \ref{corolla}}\label{sec:3}
\noindent Let $W(P)$ be the Wronskian $W(v_1,\ldots,v_N)(P)$ of the
basis $\{v_i\}_{i\in I_N}$ at a generic point $P\in C$, and
$\hat{W}_{ij}(P):=W(v_1,\ldots,v_{i-1},v_j,v_{i+1},\ldots,v_N)(P)$.

\begin{lemma}\label{illemma} The $(g-2)(g-3)/2$ linearly independent relations \begin{equation}\label{lemma}v_i(z)
W(P) =\sum_{j=1}^Nv_j(z) \hat{W}_{ji}(P)\ ,\end{equation}
$i=N+1,\ldots,M$, hold $\forall z\in C$.\end{lemma}

\begin{proof} Immediate consequence of the Cramer rule. \end{proof}

\begin{remark} The ratio $\hat{W}_{ij}(P)/W(P)$ does not depend on $P$.\end{remark}

\begin{remark} Since for $i\in I_g$ \begin{equation}\label{devv}\begin{cases}v_j(p_i)=\delta_{ji}\ , &j\in I_g\ ,\\
&\\ v_j(p_i)=0\ ,&j=g+1,\ldots,M\ ,\end{cases}\end{equation} it
follows that for $z=p_i$ Eq.\eqref{lemma} gives
$\hat{W}_{ij}(P)=0$ for $i\in I_g$ and
$j=N+1,\ldots,M$.\end{remark}

\begin{remark}\label{lastrem} The Wronskians in the expansion \eqref{lemma} can be replaced by the corresponding determinant $\det v_j(x_i)$,
where $x_1,\ldots,x_{3g-3}$ are arbitrary points on $C$.
\end{remark}

\begin{nproof}{of Theorem \ref{teo}} Assume that
$\det\,\omega(p_1,\ldots,p_g)\ne 0$. Define
$$x_i:=\begin{cases}
p_i\ ,&i\in I_g\ ,\\ &\\ q_{i-g}\ ,&i=g+1,\ldots,N+1\
.\end{cases}$$ Fix $k,l$ with $3\le k<l\le g$ and consider the
matrix

$$\begin{pmatrix}
v_1(x_1) &\ldots &v_{N}(x_1) &\sigma_k(x_1)\sigma_l(x_1)\\ \vdots
&\ddots &\vdots &\vdots\\ v_1(x_{N+1}) &\ldots &v_{N}(x_{N+1})
&\sigma_k(x_{N+1})\sigma_l(x_{N+1})\end{pmatrix}\ .$$

\noindent By \eqref{devv}, this matrix has ${\rm
diag}\,(1,\ldots,1)$ in the $g\times g$ upper left corner, $0$ in
the $g \times (2g-2)$ upper right corner and
$(\det\,\omega(p_1,\ldots,p_g))^{-2}A(k,l)$ in the $(2g-2)\times
(2g-2)$ lower right corner. On the other hand, by Lemma
\ref{illemma} and by Remark \ref{lastrem} the determinant of this
matrix vanishes. Since such relations hold for $(p_1,\ldots,p_g)$ in
a dense subset of $C^g$, they hold $\forall(p_1,\ldots,p_g)\in C^g$
and the theorem follows.\qed\end{nproof}

\begin{nproof}{of Corollary \ref{corolla}} Divide the relations in Theorem \ref{teo} by $\Delta_{i\, 2g-2}$ and
note that
$$a_{mn,r}=\sum_{i,j=1}^g D_{mi}D_{nj}\omega_i(q_r)\omega_j(q_r)\ .$$ Independence of the points
$q_i$, $i\in I_{2g-2}$, follows by noting that the coefficients of
$\omega_i\omega_j$ in the relations are functions of $q_i$ with no
zeroes or poles.\qed\end{nproof}

Define
$$(\1_i,\2_i):=\begin{cases}(i,i)\ , &1\le i\le g\ ,\\
(1,i-g+1)\ , &g+1\le i\le 2g-1\ ,\\
(2,i-2g+3)\ , &2g\le i\le 3g-3\
,\\ \hfill\vdots \hfill& \hfill\vdots\hfill\\
(g-1,g)\ , &i=g(g+1)/2\ ,\end{cases}$$ so that $\1_i\2_i$ is the
$i$-th element in the $M$-tuple
$(11,22,\ldots,gg,12,\ldots,1g,23,\ldots)$. $\forall u\in\CC^g$ and
for all the $g\times g$ matrices $A$, set
\begin{equation}\label{doubleindex}uu_i:=u_{\1_i}u_{\2_i}\ ,\quad\;
(AA)_{ij}:=\frac{A_{\1_i\1_j}A_{\2_i\2_j}+A_{\1_i\2_j}A_{\2_i\1_j}}{1+\delta_{\1_j\2_j}}\ ,\quad\;
A_i:=A_{\1_i\2_i}\ ,
\end{equation}
$i,j\in I_M$. In the following we will repeatedly make use of the identities
\begin{equation}\label{AA}\sum_{i,j=1}^gf(i,j)=\sum_{i\le j}^g\frac{f(i,j)+f(j,i)}{1+\delta_{ij}}
=\sum_{k=1}^M\frac{f(\1_k,\2_k)+f(\2_k,\1_k)}{1+\delta_{\1_k\2_k}}\
.\end{equation} In particular, if $f(i,j)=f(j,i)$, then
\begin{equation}\label{BB}\sum_{i,j=1}^gf(i,j)=\sum_{k=1}^M(2-\delta_{\1_k\2_k})f(\1_k,\2_k)\ ,\end{equation}
where we used the identity
$$2-\delta_{ij}=\frac{2}{1+\delta_{ij}}\ .$$
With this notation, and observing that
$\sigma_i=\sum_{j=1}^g[\omega]^{-1}_{ij}\omega_j$, we have
\begin{equation}\label{vomega}v_i=\sigma\sigma_i=
\sum_{j=1}^M([\omega]^{-1}[\omega]^{-1})_{ij}\omega\omega_j\ ,
\qquad i\in I_M\ .\end{equation} Set $w_{ij}:=W_{ij}/W$, where
${W}_{ij}(P):=W(v_1,\ldots,v_{i-1},\omega\omega_j,v_{i+1},\ldots,v_N)(P)$,
 and note that
\begin{equation}\label{teor}\omega\omega_i=\sum_{j=1}^N
w_{ji}v_j\ ,\qquad i\in I_M\ .\end{equation}

\section{Siegel's induced measure on $\hat{\mathcal M}_g$ and Bergman reproducing kernel}\label{sec:4}

Let $$\Sieg_g:=\{Z\in M_g(\CC)\,|\, {}^tZ=Z,\ Y>0\}\ , \qquad Y:={\Im}Z\ ,$$ be
the Siegel upper half-space and $\{\alpha_1,\ldots,\alpha_g,\beta_1,\ldots,\beta_g\}$ a
symplectic basis of $H_1(C,\ZZ)$. Denote by
$\{\omega_i\}_{i\in I_g}$ the basis of $H^0(K_C)$, dual of
$H_1(C,\ZZ)$, so that $\oint_{\alpha_i}\omega_j=\delta_{ij}$,
$i,j\in I_g$. Let $\tau_{ij}:=\oint_{\beta_i}\omega_j\in \Sieg_g$ be the Riemann
period matrix of $C$. Under the symplectic transformation
$$\begin{pmatrix}\tilde\alpha \\ \tilde\beta\end{pmatrix}=\begin{pmatrix}D & C \\ B & A\end{pmatrix}
\begin{pmatrix}\alpha\cr\beta\end{pmatrix}\ ,\qquad\qquad \begin{pmatrix}A & B\\ C & D\end{pmatrix}\in Sp(2g,\ZZ)\ ,$$
we have $\tilde\omega={}^t(C\tau+D)\cdot\omega$, with
$\tilde\tau_{ij}$ and $\tau_{ij}$ related by
the modular transformation
$$
\tilde\tau=(A\tau+B)\cdot(C\tau+D)^{-1}\ .
$$
Note that the basis $\{\sigma_i\}_{i\in I_g}$, defined in
Eq.\eqref{newbasis}, is independent of the choice of the basis of
$H_0(K_C)$ and therefore is modular invariant.

The Siegel metric
\begin{equation}\label{siegelon}
ds^2:=\Tr\, (Y^{-1}dZ Y^{-1}d\bar Z)\ ,\end{equation} defines the volume form
$$d\nu=\frac{i^M}{2^g}\frac{ \bigwedge_{i\le j}^g ({\rm d} Z_{ij}\wedge{\rm d}\bar Z_{ij})  }{
\det Y^{g+1} } \ . $$

We use the indexing introduced in Eq.\eqref{doubleindex} to express the Siegel metric on $\Sieg_g$
where now the matrix elements $Z_{ij}$, $i,j\in I_M$, are seen as the components of the $M$-dimensional vectors
$Z:=(Z_1,\ldots,Z_M)$.

\begin{proposition}
\begin{equation}\label{sieggmetr}ds^2=\sum_{i,j=1}^Mg_{ij}dZ_id\bar Z_j\ ,\end{equation}
where
\begin{equation}\label{metric}
g_{ij}:=(2-\delta_{\1_i\2_i})(Y^{-1}Y^{-1})_{ij}\ .
\end{equation}
\end{proposition}

\begin{proof} By \eqref{AA} and \eqref{BB}
\begin{align*}ds^2&=\sum_{i,j,k,l=1}^gY^{-1}_{ij}d Z_{jk}Y^{-1}_{kl}d\bar Z_{li}\\
&=\sum_{i,l=1}^gd\bar Z_{li}\sum_{m=1}^M\frac{Y^{-1}_{i\1_m}Y^{-1}_{l\2_m}+
Y^{-1}_{i\2_m}Y^{-1}_{l\1_m}}{1+\delta_{\1_m\2_m}}d Z_{\1_m\2_m}\\
&=\sum_{m,n=1}^M(2-\delta_{\1_n\2_n})d\bar Z_{\1_n\2_n}\frac{Y^{-1}_{\1_n\1_m}
Y^{-1}_{\2_n\2_m}+
Y^{-1}_{\1_n\2_m}Y^{-1}_{\2_n\1_m}}{1+\delta_{\1_m\2_m}}d Z_{\1_m\2_m}\\
&=\sum_{m,n=1}^M(2-\delta_{\1_n\2_n})(Y^{-1}Y^{-1})_{nm}d Z_md\bar Z_n\ .
\end{align*}\end{proof}

Let $k$ be the Kodaira-Spencer map identifying the
quadratic differentials on $C$ with the fiber of the cotangent of
the Teichm\"uller space at $C$. We have
$$k(\omega_i\omega_j)={(2\pi i)}^{-1}d\tau_{ij} \ . $$ By Corollary
\ref{corolla} it follows that
\begin{equation}\sum_{i,j=1}^g\frac{\det\,A_{ij,r}(k,l)}{ \Delta_{r\, 2g-2}}d\tau_{ij}=0\ ,
\label{hnice}\end{equation}
$3\le k<l\le g$. Set $d\tau_i:=d\tau_{\1_i\2_i}$, $i\in I_M$.
Eq.\eqref{metric} yields an explicit expression for the volume form on
$\hat\M_g\hookrightarrow\Sieg_g/Sp(2g,\ZZ)$ induced by the modular invariant Siegel metric on
$\Sieg_g$. Set $\tau_2:=\Im\tau$, and let $\left|\tau_2^{-1}\tau_2^{-1}\right|^{i_1\ldots
i_N}_{j_1\ldots j_N}$, with $i_k,j_k$, $k\in I_N$, distinct elements of $I_M$, be the determinant of the $N\times N$ submatrix of
$(\tau_2^{-1}\tau_2^{-1})_{ij}$, built by taking the rows $i_1,\ldots,i_N$ and the columns
$j_1,\ldots,j_N$.

\begin{theorem} The volume form on $\hat\M_g$ induced by the Siegel metric is
\begin{equation}\label{basicc}d\nu_{|\hat\M_g}=\Bigl(\frac{i}{2}\Bigr)^N\sum_{\substack{i_N>\ldots>i_1=1\\
j_N>\ldots>j_1=1}}^M\left|\tau_2^{-1}\tau_2^{-1}\right|^{i_1\ldots
i_N}_{j_1\ldots j_N}\prod_{k=1}^N(2-\delta_{\1_{i_k}\2_{i_k}})\bigwedge_1^N(d\tau_{i_k}\wedge
d\bar\tau_{j_k})\ ,
\end{equation}
so that
\begin{equation}\label{basicccc}
{\rm Vol}(\hat\M_g)=\int_{\hat\M_g}d\nu_{|\hat\M_g}\ .
\end{equation}
\end{theorem}

\begin{proof}
Let $$\omega:=\frac{i}{2}\sum_{i,j=1}^Mg_{ij}dZ_i\wedge d\bar Z_j\ ,$$
be the $(1,1)$-form associated to the Siegel metric on $\Sieg_g$. By Wirtinger's theorem \cite{GH}, the
volume form on a $d$-dimensional complex submanifold $S$ is
$$\frac{1}{d!}\omega^d\ ,$$
so that the volume of $S$ is expressed as the integral over $S$ of a globally defined differential form on $\Sieg_g$.
Set $g^\tau_{ij}:=(2-\delta_{\1_i\2_i})(\tau_2^{-1}\tau_2^{-1})_{ij}$, $i,j\in I_M$, and note that
\begin{align*}d\nu_{|\hat\M_g}=&\frac{i^N}{2^NN!}\sum_{\substack{i_1,\ldots,i_N=1\\j_1,\ldots,j_N=1}}^M
\prod_{k=1}^Ng^\tau_{i_kj_k}\bigwedge_{k=1}^N(d\tau_{i_k}\wedge
d\bar \tau_{j_k})\\
=&\frac{i^N}{2^NN!}\sum_{\substack{i_N<\ldots<i_1=1\\j_N<\ldots<j_1=1}}^M\sum_{r,s\in\perm_N}
\epsilon(r)\epsilon(s)\prod_{k=1}^Ng^\tau_{i_{r(k)}j_{s(k)}}\bigwedge_{k=1}^N(d\tau_{i_k}\wedge
d\bar \tau_{j_k})\ ,
\end{align*}
where $\perm_N$ is the group of permutations of $N$ elements and $\epsilon(s)$ is the sign of the
permutation $s$. The theorem then follows by the identity
$$\sum_{r,s\in\perm_N}
\epsilon(r)\epsilon(s)\prod_{k=1}^Ng^\tau_{i_{r(k)}j_{s(k)}}=N!\left|\tau_2^{-1}\tau_2^{-1}\right|^{i_1\ldots
i_N}_{j_1\ldots j_N}\prod_{k=1}^N(2-\delta_{\1_{i_k}\2_{i_k}})\ .$$
\end{proof}

Petri's basis of $H^2(K_C^2)$ corresponds, through the Kodaira-Spencer map, to a basis for the cotangent space of the
Teichm\"uller space.
 Setting $d\,\Xi_i:=2\pi i\,k(v_i)$, $i\in I_N$,
it follows by Eq.\eqref{teor} that
\begin{equation}\label{domega}d\tau_i=\sum_{j=1}^Nw_{ji}
d\,\Xi_j \ ,\qquad i\in I_M\ .\end{equation}

\begin{corollary} Fix the points $p_1,\ldots,p_g\in C$ in general position,
so that $\{v_i\}_{i\in I_N}$, is a basis of $H^0(K_C^2)$. The metric on $\hat\M_g$ induced by the Siegel metric
is \begin{equation}\label{ds} ds_{|\hat\M_g}^2=
\sum_{i,j=1}^Ng^{\Xi}_{ij}d\,\Xi_id\,\bar\Xi_j \ ,\end{equation} where
$g^{\Xi}_{ij}:=\sum_{k,l=1}^M(2-\delta_{\1_k\2_k})w_{ik}(\tau_2^{-1}\tau_2^{-1})_{kl}
\bar{{w}}_{jl}$. \end{corollary}

\begin{proof} Immediate. \end{proof}

By using a suitable basis of $H^0(K_C^2)$ and its image
under the Kodaira-Spencer map, it turns out that the metric $g$ is
related to the Bergman reproducing kernel. Fix the points $z_1,\ldots,z_N\in C$
satisfying the conditions of Proposition \ref{prop:base}. The basis
$\{\gamma_i\}_{i\in I_N}$ of $H^0(K_C^2)$, with $\gamma_i\equiv\gamma^2_i$, $i\in I_N$, defined by
Eq.\eqref{basendiff} in the case $n=2$, satisfies the
relations
$$\omega\omega_i=\sum_{j=1}^N\omega\omega_i(z_j)\gamma_j\ ,\qquad
v_i=\sum_{j=1}^Nv_i(z_j)\gamma_j\ ,\quad i\in I_M\ .$$ Set
$\Gamma_i:=(2\pi i)^{-1}k(\gamma_i)$ and $[v]_{ij}:=v_i(z_j)$, $i,j\in I_N$.

\begin{theorem}
\begin{equation}\label{siegelbergman}
ds^2_{|\hat\M_g}= \sum_{i,j=1}^NB^2(z_i,\bar
z_j)d\,\Gamma_id\,\bar\Gamma_j\ ,
\end{equation}
in particular, the Siegel induced modular invariant volume form on
$\hat\M_g$ is
\begin{equation}\label{volumee}d\nu_{|\hat\M_g}=\Bigl(\frac{i}{2}\Bigr)^N\det\,
B^2(z_i,\bar z_j)\,\textstyle{\bigwedge_1^N}(d\,\Gamma_i\wedge d\,\bar \Gamma_i) \
,
\end{equation}
where
$$
B(z,\bar w):=\sum_{i,j=1}^g\omega_i(z)(\tau_2^{-1})_{ij}\bar\omega_j(w)\ ,$$
$z,w\in C$, is the Bergman reproducing kernel. \end{theorem}

\begin{proof} Use $d\tau_i=\sum_{j=1}^N\omega\omega_i(z_j)d\Gamma_j$, $i\in I_g$,
and the identity
$$\sum_{k,l=1}^M(2-\delta_{\1_k\2_k})\omega\omega_k(z_i)(\tau_2^{-1}\tau_2^{-1})_{kl}
\bar\omega\bar\omega_{l}(z_j)=B^2(z_i,\bar
z_j)\ ,\qquad i,j\in I_N\ . $$
Note that by \eqref{ds} $\sum_{k,l=1}^N[v]_{ki}g^{\Xi}_{kl}[\bar v]_{lj}=B^2(z_i,\bar z_j)$,
which also follows by \eqref{teor}. \end{proof}

\section{Fay's trisecant identity from the distinguished basis of $H^0(K_C^n)$}\label{sec:5}

Let $I_i(p):=\int_{p_0}^p\omega_i$,
$p_0,p\in C$, $i\in I_g$, be the Abel-Jacobi map, which extends to a map from
divisors of $C$ to the Jacobian
$J(C):=\CC^g/(\ZZ_g+\tau\ZZ^g)$. We consider
Riemann $\theta$-functions $\theta(D+e):=\theta (I(D)+e,\tau)$,
$e\in J(C)$, evaluated at some $0$-degree divisor $D$ of $C$. By the
Riemann vanishing theorem, there is a divisor class $\Delta$ of
degree $g-1$ with $2\Delta = K$, such that $-I(\Delta)$ is the
vector of Riemann constants. Let $E(z,w)$ be the prime
form, and set $$\sigma(z):=\exp\bigg(-\sum_{i=1}^g\oint_{\alpha_i}\omega_i(w)\ln E(z,w)\bigg)\ .$$

\begin{proposition}\label{thdettheta} Fix $n\in\NN$ and let
$\{\phi_i^n\}_{i\in I_{N_n}}$ be an arbitrary basis of
$H^0(K_C^n)$, $n\geq 1$.
Let $y,x_1,\ldots,x_{N_n}$ be
arbitrary points of $C$. Then, for $n=1$
\begin{equation}\label{dettheta}\det\phi^1_i(x_j)=\kappa_1[\phi^1]\frac{\theta\bigl(\sum_{1}^gx_i-y-\Delta\bigr)
\prod_{i<j}^gE(x_i,x_j)\prod_1^g\sigma(x_k)}{
\sigma(y)\prod_1^gE(y,x_i)}\ ,\end{equation} whereas for $n>1$
\begin{equation}\label{detthetaii}\det\phi^{n}_i(x_j)=\kappa_{n}[\phi^n]
\theta\Bigl(\sum_{1}^{N_n} x_i-(2n-1)\Delta\Bigr)\prod_{i<j}^{N_n}
E(x_i,x_j)\prod_{i=1}^{N_n}\sigma(x_i)^{2n-1}\ ,\end{equation} where
$\kappa_1[\phi^1]$ and $\kappa_n[\phi^n]$ are constants depending only
on the choice of the bases. \end{proposition}

\begin{proof} $\kappa_1[\phi^1]$ is a nowhere vanishing
section in $x_j$, $j\in I_g$, and
$\theta(\sum_{1}^gx_i-y-\Delta)=0$ for $y=x_1,\ldots,x_g$, so that
it is also a nowhere vanishing section in $y$ and
since it has trivial monodromy it must be a constant.
Eq.\eqref{detthetaii} follows by a similar proof.
\end{proof}

Set $w:=\sum_1^{N_n}p_i-(2n-1)\Delta$, $n>1$, and assume that $p_1,\ldots,p_{N_n}\in C$
satisfy the hypothesis of Proposition \ref{prop:base}.
By Proposition \ref{thdettheta} we have
\begin{equation}\label{glindiff}\gamma^n_i(z)=\frac{\theta(w+z-p_i)\sigma(z)^{2n-1}\prod_{^{k=1}_{k\neq
i}}^{N_n}E(z,p_k)}{ \theta(w)\sigma(p_i)^{2n-1}\prod_{^{k=1}_{k\neq
i}}^{N_n}E(p_i,p_k)}\ , \qquad i\in I_{N_n}\ .\end{equation}

\begin{theorem}\label{Faytris} Propositions \ref{prop:base} and \ref{thdettheta} imply
the Fay trisecant identity \cite{Fay}
$$\frac{\theta(w+\sum_{i=1}^m(x_i-y_i))\prod_{i<j}E(x_i,x_j)E(y_i,y_j)}{
\theta(w)\prod_{i,j}E(x_i,y_j)}=(-)^{\frac{m(m-1)}{2}}\det\nolimits_{ij}\frac{\theta(w+x_i-y_j)}{
\theta(w)E(x_i,y_j)}\ ,$$ $m\ge 2$, $\forall
x_1,\ldots,x_m,y_1,\ldots,y_m\in C$, $w\in J(C)$.
\end{theorem}

\begin{proof}
Fix $m\ge 2$, $x_1,\ldots,x_m,y_1,\ldots,y_m\in C$ and
$w\in J(C)$, with $\theta(w)\neq 0$. Choose
$y_1,\ldots,y_m$ distinct, otherwise the identity
is trivial. Set $p_i:= y_i$, $i\in I_m$, and fix $n\in \ZZ$,
$N_n\ge m$, and $p_{m+1},\ldots,p_{N_n}\in C$, so that
$w=\sum_1^{N_n}p_i-(2n-1)\Delta$. By Jacobi inversion theorem
such a choice is always possible, provided that $N_n-m\ge g$. It is also clear that, for $n$ large enough,
$p_{m+1},\ldots,p_{N_n}$ can be chosen pairwise
distinct and distinct from $y_1,\ldots,y_m$. Eq.\eqref{detthetaii}
implies that, by construction, $\det\phi_i^n(p_j)\neq 0$, for any
basis $\{\phi_i^n\}_{i\in I_{N_n}}$ of $H^0(K_C^n)$,
since the points $p_1,\ldots,p_{N_n}$ are pairwise distinct and
$\theta(w)\neq 0$. Therefore, one can define a basis
$\{\gamma^n_i\}_{i\in I_{N_n}}$ of $H^0(K_C^n)$ by \eqref{basendiff}
and consider $\det\gamma^n(x_1,\ldots,x_m,p_{m+1},\ldots,p_{N_n})$, which can be
expressed by \eqref{glindiff} or
 by \eqref{detthetaii}, with
$\kappa_n[\gamma^n]$ determined by applying \eqref{detthetaii} to
 $\det\gamma^n_i(p_j)=1$. Comparing these formulas the
theorem follows.
\end{proof}

%

\bibliographystyle{amsplain}

\end{document}